\documentclass[11pt]{article}
\usepackage{amssymb}
\usepackage{enumerate}
\usepackage{enumitem}
\usepackage{mathrsfs}
\usepackage{latexsym,amsmath,amssymb,amsfonts,epsfig,graphicx,cite,psfrag}
\usepackage{eepic,color,colordvi,amscd}

\newcommand{\qed}{\hfill $\Box $}
\newcommand{\pf}{\noindent {\bf Proof.} }
\newtheorem{theorem}{Theorem}
\newtheorem{lemma}[theorem]{Lemma}

\begin{document}	
	\title{Anti-Ramsey Number of Stars in  3-uniform hypergraphs\thanks{Supported by  National Key Research and Development Program of China 2023YFA1010203 and the National Natural
Science Foundation of China under grant No.12271425}}
	\author{Hongliang Lu,\ Xinyue Luo and Xinxin Ma\\School of Mathematics and Statistics\\
Xi'an Jiaotong University\\
Xi'an, Shaanxi 710049, China\\
\smallskip\\
}
\date{}

\maketitle

\date{}

\maketitle
	
	\begin{abstract}
		An edge-colored hypergraph is called \emph{a rainbow hypergraph} if all the colors on its edges are distinct. Given two  positive integers $n,r$ and an $r$-uniform hypergraph $\mathcal{G}$, the anti-Ramsey number $ar_r(n,\mathcal{G})$ is defined to be the minimum number of colors $t$ such that there exists a rainbow copy of $\mathcal{G}$ in any exactly $t$-edge-coloring of the complete $r$-uniform hypergraph of order $n$. Let $ \mathcal{F}_k $ denote the 3-graph ($k$-star) consisting of $k$ edges sharing exactly one vertex. Tang, Li and Yan \cite{YTG} determined the value of $ar_3(n,\mathcal{F}_3)$ when $n\geq 20$. In this paper, we determine the anti-Ramsey number $ar_3(n,\mathcal{F}_{k+1})$, where $k\geq 3$ and $n> \frac{5}{2}k^3+\frac{15}{2}k^2+26k-3$.

	\end{abstract}
	\begin{flushleft}
\end{flushleft}
	
\section{Introduction}
For a set $S$ and a positive integer $k$, we use ${S\choose k}$ to denote the collection of all possible subsets of $k$ elements of $S$. A \textit{hypergraph} $ \mathcal{F}=(V(\mathcal{F}), E(\mathcal{F})) $ consists of a vertex set $ V(\mathcal{F}) $ and an edge set $ E(\mathcal{F}) $, where each edge in $ E(\mathcal{F}) $ is a non-empty subset of $ V(\mathcal{F}) $.  The number of edges of $\mathcal{F}$ is  denoted by $e(\mathcal{F})$, that is, $e(\mathcal{F}):=|E(\mathcal{F})|$.  
If $ |e|=r $ for any $e\in E(\mathcal{F})$, then $\mathcal{F}$ is called an \textit{$ r$-uniform   hypergraph}  (or $r$-graph,  for simplicity). 
For $u\in V(\mathcal{F})$, let $N_\mathcal{F}(u):=\{e\ |\ e\subseteq V(\mathcal{F})\setminus\{u\}\ \mbox{and}\ e\cup \{u\}\in E(\mathcal{F})\}$ be the neighborhood of $u$ in $\mathcal{F}$. The \textit{degree} of $u$ in $\mathcal{F}$, denoted by $d_\mathcal{F}(u)$, is the size of $N_{\mathcal{F}}(u)$.
 For $X\subseteq V(\mathcal{F})$, we define $\mathcal{F}-X$ as the subhypergraph of $\mathcal{F}$ obtained by removing all vertices in $X$ and all  edges intersecting with $X$ in $\mathcal{F}$. 
 Similarly, if $Y\subseteq E(\mathcal{F})$,  we use 
$\mathcal{F}-Y$ to denote the hypergraph resulting from deleting all the edges in  $Y$ from $\mathcal{F}$.  
When $X=\{x\}$, $Y=\{e\}$, we respectively write $\mathcal{F}-X=\mathcal{F}-x$, $\mathcal{F}-Y=\mathcal{F}-e$. Specifically, to avoid confusion, for an edge  $e\in E(\mathcal{F})$, we use $\mathcal{F}-V(e)$ to denote the subgraph of $\mathcal{F}$ obtained  by removing all vertices in $e$ and all  edges intersecting with $e$ in $\mathcal{F}$.
For a non-empty subset $X\subseteq V(\mathcal{F})$, let $\mathcal{F}[X]$ denote the subgraph \emph{induced} by $X$. For two disjoint sets $U$ and $W$, we use $U\times W$ to denote the collection of  2-sets that  intersect $U$ and $W$. That is, $U\times W=\{\{x,y\}\ |\ x\in U\mbox{ and }y\in W\}$.


 	For two vertex-disjoint hypergraphs $\mathcal{F}, \mathcal{H}$, the \textit{union} of $\mathcal{F}$ and $\mathcal{H}$ denoted by $ \mathcal{F}\cup \mathcal{H}$   is the hypergraph with vertex set $V(\mathcal{F})\cup V(\mathcal{H})$ and edge set $E(\mathcal{F})\cup E(\mathcal{H})$. When there is no confusion, for $T\subseteq V(\mathcal{F})$, we also use $N_{\mathcal{F}}(T)$   to denote the $(k-|T|)$-graph with vertex set $V(\mathcal{F})-T$ and edge set $N_{\mathcal{F}}(T)$. A \emph{matching} in a hypergraph $H$ is a set of pairwise disjoint edges in $H$, and we use $\nu(H)$ to denote the maximum size of a matching in $H$.
 	
 	A $t$-\textit{edge-coloring} of a hypergraph is an assignment of $t$ colors to its edges, and an \textit{exactly $t$-edge-coloring} uses all $t$ colors. An edge-colored graph is called \textit{rainbow} if all edges have distinct colors. Let $ [n]=\{1,...,n\}$.  The complete $ r $-graph with order $n$ is denoted by $ K_n^r$ . Given a positive integer $ n $ and a hypergraph $ \mathcal{F} $, the \textit{anti-Ramsey number} $ ar_r(n,\mathcal{F}) $ is the minimum number of colors $t$ such that each edge-coloring of $K_n^r$ with exactly $t$ colors contains a rainbow copy of $\mathcal{F}$.  Given an edge-coloring $C$ of $\mathcal{F}$, the colored  hypergraph  is \emph{$F$-free} if $\mathcal{F}$ has no rainbow subhypergraph which is isomorphic to $F$.

Given a hypergraph   $H$ and a family of hypergraphs $\mathcal{H}$,   $H$ is called \emph{$\mathcal{H}$-free} if for any $F\in \mathcal{H}$, $H$ does not contain  $ F$ as a subhypergraph. 
The Tur\'an number of a fmaily of $r$-graphs $\mathcal{F}$, written $ex_r(n, \mathcal{F})$, is the largest possible number of edges in an $\mathcal{F}$-free $r$-graph on $n$ vertices. When $\mathcal{F}=\{F\}$, we use $ex_r(n, F)$ instead of $ex_r(n, \{F\})$.

Anti-Ramsey numbers were introduced by Erd\H{o}s, Simonovits and S\'{o}s \cite{PMV} in 1973. They found that the anti-Ramsey numbers are closely related to the Tur\'{a}n numbers.  For an  $r$-graph $F$, there is a natural lower bound of $ar_r(n,F)$ in terms of  Tur\'{a}n number as follow,
\begin{align}{\label{rainbow}}
	ar_r(n,F)\geq ex_r(n,\{F-e:e\in E(F)\})+2.
\end{align}
This trivial lower bound is easily obtained by coloring a rainbow Tur\'{a}n extremal $r$-graph for $\{F-e:e\in F\}$ in $K_n^r$, and the remaining edges with an additional color \cite{TGG}. 

In 1973, Erdős, Simonovits, and Sós \cite{PMV} proved that there exists an integer $n_0(p)$ such that for all $n > n_0(p)$, the anti-Ramsey number $ar(K_p,n)$ satisfies the equation $ar(K_p,n)=ex(n,K_{p - 1})+2$. Later, Montellano-Ballesteros and Neumann-Lara \cite{MN}, as well as Schiermeyer \cite{S}, independently extended this result to cover all values of $n$ and $p$ where $n>p\geq3$. 
Jiang \cite{Jiang} and Montellano-Ballesteros \cite{MB} independently determined the anti-Ramsey numbers for stars in graphs. A variety of results regarding the anti-Ramsey numbers of 2-graphs have been achieved. These results cover various graph structures such as paths, cycles, and matchings.
 For a comprehensive overview, we recommend the survey paper \cite{SCK}. Regarding  3-graphs, Guo, Lu, and Peng \cite{MHX}  established the exact value of the anti-Ramsey number for matchings. Gu, Li, and Shi \cite{RJY} investigated the anti-Ramsey numbers of paths and cycles in hypergraphs. For other related results on the anti-Ramsey numbers of paths, cycles, and matchings in hypergraphs, interested readers are referred to \cite{PA,TGG,YT,YTG,MB,LS25}. 
 In this paper, we aim to investigate the anti-Ramsey numbers of stars of 3-graphs.

Let $ \mathcal{F}_k $ ($k$-star) denote the 3-graph  consisting of $k$ edges sharing exactly one vertex, called the core of the star. 
Let $ f(n,k) $ denote the maximum number of edges in $ 3 $-graph without  $ k $-stars. $f(n,2)$ was determined exactly by Erd\H{o}s and S\'os\cite{VT}. 
Duke and Erd\H{o}s \cite{DE} established linear lower and upper bounds on \(f(n,k)\) for fixed \(k\) (where the bounds scale linearly with \(n\)). These bounds were subsequently improved in \cite{Fr83}, where the exact value of \(f(n,3)\) was determined for all integers \(n\geq 54\). Further refinements were made in \cite{Chung}, yielding bounds that are nearly best possible. Finally, Chung and Frankl \cite{CF} derived the exact value of \(f(n,k)\) for 3-graphs in the regime \(n\geq \frac{5}{2}k^3\).




\begin{theorem}[Erd\H{o}s and S\'os, \cite{VT}]
	For all $n\geq3$, 
	$$f(n,2)=\begin{cases}
		n, & \mbox{if}\ n\equiv0\pmod 4, \\
		n-1, &\mbox{if}\ n\equiv1\pmod 4, \\
		n-2, & \mbox{if}\ n\equiv2,3\pmod 4.\\		
		\end{cases}$$
	\end{theorem}

 \begin{theorem}[Chung and Frankl, \cite{CF}]{\label{star}}
 	Suppose that $k\geq3$ is odd and $n>k(k-1)(5k+2)/2$, then $ f(n,k)=(n-2k)k(k-1))+2{k\choose 3} $. Moreover, a $3$-graph $\mathcal{F}$ has $f(n,k)$ edges and contains no $k$-star if and only if $\mathcal{F}$ is isomorphic to $\mathscr{F}_k^o$.
  \end{theorem}

  \begin{theorem}[Chung and Frankl, \cite{CF}]{\label{star}}
 	Suppose that $k\geq 4$ is even and $n>2k^3-9k+7$, then $ f(n,k)=\frac{1}{2}nk(2k-3)-\frac{1}{2}(2k^3-9k+6)$. Moreover, a $3$-graph $\mathcal{F}$ has $f(n,k)$ edges and contains no $k$-star if and only if $\mathcal{F}$ is isomorphic to $\mathscr{F}_k^e$.
  \end{theorem}
  Chung and Frankl \cite{CF} proposed the following extremal graph construction. 
  
\noindent\textit{The construction of $\mathscr{F}_k^o$ }: Let $k$ be odd and let $ S $ and $R$ be two disjoint sets of $ [n] $ of size $k$. Consider the 3-graph $\mathscr{F}_k^o$  with vertex set $[n]$ and edge set 
\[
E(\mathscr{F}_k^o)=\{T\in \binom{V}{3}:|T\cap S|\geq 2 \ \mbox{ and }\ |T\cap R|= \emptyset\} \cup \{T\in \binom{V}{3}:|T\cap R|\geq 2 \ \mbox{ and }\ |T\cap S|= \emptyset\} .
\]
\textit{The construction of $\mathscr{F}_k^e$}: Let $k$ be even. Let $G_k$ be 
the 2-graph with $2k-1$ vertices, denoted by $x_1,\ldots,x_{k-1}, y_1,\ldots,y_{k-1}$ and $z$. The edge set of $G_k$ consists of all the pairs $(x_i, y_j)$, except for $(x_i, y_i)$
with $2i>k$ together with the pairs $(x_i, z), (y_i, z)$ with $2i > k$. It is easy to see that   it has all degrees equal $k-1$
except for the degree of $z$, which is $k-2$. Let $\mathscr{F}_k^e$ denote the 3-graph on $n$ vertices  such that each edge either intersects the vertex set $V(G_k)$ in an edge of $G_k$ or contains two distinct edges of $G_k$, together with all the triples of the form $\{x_i,y_i,z\}$ for $1\leq i\leq k/2$.

The anti-Ramsey number $ar_3(n,\mathcal{F}_2)=2$ follows immediately from the definition. Tang, Li and Yan \cite{YTG} determined the value of $ar_3(n,\mathcal{F}_3)$ when $n\geq 20$. 

\begin{theorem}[Tang, Li and Yan \cite{YTG}]\label{tang}
	For all $n\geq 20$, 
$$ar_3(n,\mathcal{F}_3)=\left\{
\begin{array}
{ll}f(n,2)+2, & \mbox{if}\ n\equiv0\pmod 4, \\
f(n,2)+2, & \mbox{if}\ n\equiv1\pmod 4, \\
f(n,2)+3, & \mbox{if}\ n\equiv2,3\pmod 4.
\end{array}\right.$$

\end{theorem}

In this paper, we determine the anti-Ramsey number of a star in 3-graphs for sufficiently large $n$. 
\begin{theorem}{\label{main}}
	For $k\geq 3$ and $n> \frac{5}{2}k^3+\frac{15}{2}k^2+26k-3$, $ ar_3(n,\mathcal{F}_{k+1})=f(n,k)+2$. 
\end{theorem}

\section{Technical Lemmas}

\begin{theorem}[Tutte, \cite{Tu47}]\label{tutte}
A graph $G$ has  a perfect matching if and only if 
$$o(G-S)\leq |S|\quad \mbox{for any }S\subseteq V(G),$$
where $o(G-S)$ denotes the number of connected  components of odd order in $G-S$.
\end{theorem}
A graph $G$ is called a \emph{factor-critical} graph if $G-v$ has a perfect matching for each $v\in V(G)$.
Gallai \cite{Ga63} introduced the concept of factor-critical graphs. 
By definition and Theorem \ref{tutte}, one can see that the following statement holds. 
\begin{lemma}[Gallai, \cite{Ga63}]\label{critical}
A graph $G$ of odd order is factor-critical if and only if $$o(G-S)\leq |S|\quad \mbox{for any non-emptyset  }S\subseteq V(G),$$
where $o(G-S)$ denotes the number of connected  components of odd order in $G-S$.
\end{lemma}

Here we use a ``weight function" method developed by Chung and Frankl \cite{CF} to characterize the structure of a 3-graph $\mathcal{F}$.  
Suppose that $\mathcal{F}$ is a family of 3-element subsets of the $n$-set $V=V(\mathcal{F})$. Let $P$ denote the set of all pairs of vertices in $V$. For each $\{u,v\}$ in $P$, the pair frequency is $$z(u,v)=|\{w:\{u,v,w\}\in \mathcal{F}\}|.$$
We also have 
\begin{flalign*}
	& A:=\{\{u,v\}\in P:z(u,v)\geq 2k-1\}\\
	& B:=\{\{u,v\}\in P:2k-2\geq z(u,v)\geq k\}\\
	& C:=P-A-B
\end{flalign*}
The weight function $\omega: \mathcal{F} \times P \rightarrow R$ distributes weights to pairs within each triple in $\mathcal{F}$ according to the pair frequency.

For a fixed triple $T\in\mathcal{F}$, the three pairs in $T$ are denoted by $z(p_1)\geq z(p_2)\geq z(p_3)$. The \emph{weight function} $w$ is defined as follows:
\begin{itemize}
	\item If $p_1, p_2, p_3 \in A\cup B$ or $p_1, p_2, p_3 \in B\cup C$, then $\omega(T,p_i)=\frac{1}{3}$.
	\item Suppose $p_1\in A, p_3 \in C$. If $p_2 \in A \cup B$, then $\omega(T,p_1)=\omega(T,p_2)=\frac{1}{2}, \omega(T,p_3)=0$. If $p_2 \in C$, then $\omega(T,p_1)=1, \omega(T,p_2)=\omega(T,p_3)=0$.
	\item For convenience we set also $\omega(T,p)=0$ for $p\nsubseteq T$.
\end{itemize}
Obviously, we have $$\sum_{1\leq i\leq 3}\omega(T,p_i)=1$$
and $$ \sum_{T\in \mathcal{F}}{\sum_p \omega(T,p)}=e(\mathcal{F})$$
\begin{lemma}[\cite{CF,Zhu}]\label{weight}
	For every vertex $v$ in a 3-graph $\mathcal{F}$ which does not contain a $k$-star, the following holds: $$ W_v=\sum_{p\in N_{\mathcal{F}}(v)}\omega(\{v\}\cup p,p)\leq k(k-1).$$ 
	Moreover $W_v\leq k(k-1)-\frac{2}{3}$ unless $N_{\mathcal{F}}(v)$ is the disjoint of two complete 2-graphs on $k$ vertices and every edge of $N_\mathcal{F}(v)$ is $A$-type. If $k$ is even, then one has the stronger inequality,
\[W_v\leq k(k-3/2)\]
Moreover, $W_v\leq  k(k-3/2)-1/2$ unless $N_{\mathcal{F}}(v)=K_{k-1}\cup C$, where $C$ is a factor-critical graph of order $k+1$ with degree sequence $k-1,\ldots,k-1,k-2$ or $N_{\mathcal{F}}(v)$  with maximum degree $k-1$ satisfies the following  three conditions. 
\begin{itemize}
    \item [(a)] There exists $S\subseteq V(N_{\mathcal{F}}(v))$ such that $N_{\mathcal{F}}(v)-S$ consists of isolated vertices and one factor-critical  component denoted by $F_0$ with $2k-1-2|S|\geq k+1$ vertices and degree sequence $(k-1, ..., k-1, k-2)$;
    \item [(b)]  $N_{\mathcal{F}}(v)-V(F_0)$ is the edge disjoint union of $|S|$ stars, each with maximum degree $k-1$; and
    \item [(c)] every edge of $N_{\mathcal{F}}(v)$ is in $A$, and all edges connecting $v$ to $V(N_{\mathcal{F}}(v))$ are in $C$. 
\end{itemize}
\end{lemma}
Chung and Frankl \cite{CF} first determined the value of \( f(n,k) \), though they did not characterize the corresponding extremal case: when \( N_\mathcal{F}(v) = K_{k-1} \cup C \) (where \( C \) is a factor-critical graph of order \( k+1 \) with degree sequence \( (k-1, \dots, k-1, k-2) \)). Later, Zhu et al. \cite{Zhu} adopted a similar approach to fully characterize all extremal graphs (a method analogous to the one Chung and Frankl used in \cite{CF} when handling the function \( \text{ex}_3(n,\mathcal{F}_k) \)).




Write $c(n,k):=ex_3(n,\mathcal{F}_{k})+2$. Given an edge-coloring $c:E(K^3_n)\rightarrow [c(n,k)]$ such that $c$ is surjective and the colored hypergraph denoted by $H$ contains no rainbow $\mathcal{F}_{k+1}$. For $U\subseteq V(H)$, let $Z_c(U):=\{c(e)\ |\ e\in H,U\subseteq e\}$  and let $z_c(U):=|Z_c(U)|$. When $U=\{x\}$, we denote $Z_c(U)$ and $z_c(U)$ by $Z_c(x)$ and $z_c(x)$, respectively. 
If $z_c(\{u,v\})\leq 3k$, the pair $\{u,v\}$ is \emph{good} in $H$, saying \emph{bad}, otherwise. 
	\begin{lemma}{\label{pre}}
			 There exist $2k+6$ disjoint good pairs in $H$. 
	\end{lemma}
	 
	\pf Firstly, we show that the following claim. 

    \medskip
    \textbf{Claim 1.~}For every $u$, there exist at least $n-k-1$ vertices $v\in V(H)$ such that $z_c(u,v)\leq 3k$.
   \medskip
   
By contradiction. Suppose the result does not hold. Then there exists $u\in V(H)$ and $S\in {V(H)-u\choose k+1}$ such that $z_c(u,v)\geq 3k+1$ for all $v\in S$. Write $S:=\{v_1,\ldots,v_{k+1}\}$. We choose a maximal rainbow star $\mathcal{F}_r$ in $H$ with center $u$ such that $V(\mathcal{F}_r)\cap S=r$. Since $H$ is a complete graph, we may choose $e_1$ such that $\{u,v_1\}\subseteq e_1$ and $|e_1\cap S|=1$. So we have $V(\mathcal{F}_r)\neq \emptyset$. Next we show that $r\geq k+1$, which will contradict the hypothesis. Otherwise, suppose that $r\leq k$. Then let $v\in S-V(\mathcal{F}_r)$. Since $z_c(u,v)\geq 3k+1$ and the number of edges containing $\{u,v\}$ and intersecting $V(\mathcal{F}_r)-u$ is at most $2r$, there exists an edge $f$ containing $\{u,v\}$ and $f-\{u,v\}\notin V(\mathcal{F}_r)$ and $c(f)\notin \{c(e)\ |\ e\in \mathcal{F}_r\}$. So $\{f\}\cup E(\mathcal{F}_r)$ induces a rainbow $\mathcal{F}_{r+1}$, which contradicts the choice of $\mathcal{F}_r$. This completes the proof of claim 1.

Since $n-k-1>2(2k+6)$, by Claim 1, we may greedily choose a set of $2k+6$ vertex-disjoint good pairs in $H$. This completes the proof. \qed

 The following lemma is crucial for the proof of our main theorem.

\begin{lemma}\label{degree-critical}
Let $k\geq 5$ be an integer and let $G$ be a  graph with at most $2k-1$ vertices. If  $G$ contains at most one vertex of degree $k-2$ and the rest vertices has degree $k-1$, then for any  $f\in E(G)$, $G-f$ is also factor-critical.
\end{lemma}
\pf  Otherwise, suppose that the result does not hold. Then there exists $f\in E(G)$ such that $G-f$ is not factor-critical. Then by Lemma \ref{critical}, there exists $S\subseteq V(G)$ such that $o(G-f-S)\geq |S|+1$. Recall that $\delta(G)\geq k-2$ and $G$ contains at most one vertex of degree $k-2$. 

\medskip
\textbf{Claim 1.}~$G$ is $(k-2)$-edge-connected.
\medskip

  Otherwise, suppose that there exists $M\subseteq E(G)$ such that $|M|\leq k-3$ and $G-M$ is not connected. Let $F_1$ and $F_2$ be two connected components of $G-M$ such that $|V(F_1)|\leq |V(F_2)|$. Since $\delta(G)\geq k-2$, then we have $2\leq |V(F_1)|\leq k-1$.
So we have
\begin{align*}
|M|&\geq e_{G}(V(F_1),V(G)-V(F_1))\\
&\geq |V(F_1)|(k-|V(F_1)|)-1\\
&\geq k-2,
\end{align*}
which contradicts the hypothesis that $|M|\leq k-3$. This completes the proof of Claim 1. 

By Claim 1,  we have $q\geq 2.$ Let $C_1,\ldots,C_q$ denote these odd components of $G-S-f$ such that $|V(C_1)|\geq \cdots\geq |V(C_q)|$. 

\medskip
\textbf{Claim 2.}~$|V(C_1)|\geq k$.
\medskip

Otherwise, suppose $|V(C_1)|\leq k-1$.    Then for $1\leq i\leq q$, $e_G(V(C_i),V(G)-V(C_i))\geq  k-2$ with equality if and only if the vertex of degree $k-2$ in $G$ belongs to $C_i$. Thus 
\begin{align*}
(k-1)|S|&\geq \sum_{x\in S}d_G(x)\\
&\geq  \sum_{i=1}^q e_G(S,V(C_i))\\
&\geq (k-1)(|s|+1)-3\\
&>(k-1)|S| \quad \mbox{(since $k\geq 5$)},
\end{align*}
a contradiction. This completes the proof of Claim 2. 

 By Claim 2, we have $|V(G)-V(C_1)|\leq k-1$. 

Suppose that  $\sum_{i=3}^q|V(C_i)|\geq 1$. For any $x\in V(\cup_{i=2}^q C_i)$, it follows that  $d_{G-f}(x)\leq k-3$. Recall that $G$ contains exactly one vertex of degree $k-2$. This implies that $\sum_{i=2}^q |V(C_i)|\leq 2$, which further yields $q= 3$ and $|V(C_2)|=1$. 
Then for $x\in V(C_2\cup C_3)$, 
\[
d_G(x)\leq |S|+1\leq 3,
\]
which contradicts the fact that $G$ contains exactly one vertex of degree $k-2$. 

We proceed by considering the case where $\sum_{i=3}^q|V(C_i)|=0$. This immediately implies that $q=2$ and $|S|=1$.
By Claim 1, $G$ is 2-connected. Hence $f\in E_G(V(C_1),V(C_2))$. 
If $|V(C_2)|\geq 3$, then there exist two vertices, say $x_1,x_2\in V(C_2)-V(f)$ such that $d_G(x_i)\leq k-2$, which contradicts the fact that $G$ contains exactly one vertex of degree $k-2$. We thus conclude that $|V(C_2)|=1$. Then  for $x\in V(C_2)$,
$d_G(x)\leq 2$, which once again contradicts the  condition that $\delta(G)\geq 3$. 
This completes the proof. \qed

	 \begin{lemma}\label{cycle}
	    Let $G$ be a simple graph with $|V(G)|=6$ and degree sequence $(3,3,3,3,2,2)$. Then $G$ has a Hamiltonian cycle.
	 \end{lemma}
	 
	 \pf Let $V(G)=\{w_1,w_2,w_3,w_4,u,v\}$, where $d_G(u) = d_G(v)=2$ and $d_G(w_i)=3$ for $i\in[4]$.
Without loss of generality, assume that $uw_1, vw_2\in E(G)$.

Firstly, we consider the case where $uv\in E(G)$.
The subgraph $G_1 = G-\{u, v\}$ has a degree sequence of $(3,3,2,2)$. Then $G_1$ has a path $P$ of length three that connects $w_1$ and $w_2$. The edge-set $(E(P)\cup\{uw_1,vw_2,uv\}$ forms a Hamiltonian cycle of $G$.

Next, we consider the case where $uv\notin E(G)$. If $G-\{u,v\}$ is a 4-cycle, denoted as $w_1w_2w_3w_4$, then by symmetry, we may assume that either $N_G(u)=\{w_1,w_3\}$ or $ N_G(u)=\{w_1,w_4\}$. In either of these two cases, $G$ has a Hamiltonian cycle.
Otherwise, $G-\{u,v\}$ is a not 4-cycle, which implies that $N_G(u)\cap N_G(v)\neq\emptyset$. Note that $N_G(u)\neq N_G(v)$. So, we can assume that $w_3\in N_G(u)\cap N_G(v)$. Without loss of generality, suppose that $N_G(u)=\{w_1,w_3\}$ and $N_G(v)=\{w_2,w_3\}$.  It follows that the edge set  $E(G)-\{w_4w_3,w_1w_2\}$ induces a Hamiltonian cycle of $G$.
This completes the proof. \qed

  

\section{Proof of Theorem \ref{main}}

By (\ref{rainbow}), we have $ar_3(n,\mathcal{F}_{k+1})\geq f(n,k)+2$. Therefore the lower bound is followed. Next we assume $k\geq 3$. Let $c(n,k):=f(n,k)+2$. 

For the upper bound, we prove it by contradiction. Suppose that the result does not hold. Then there exists an edge-coloring $c:E(K^3_n)\rightarrow [c(n,k)]$ such that $c$ is surjective and the colored hypergraph denoted by $\mathcal{G}$ contains no rainbow $\mathcal{F}_{k+1}$. 
By Lemma \ref{pre}, we can denote a set of $2k+6$ disjoint good pairs in $\mathcal{G}$ by $Q=\{ \{u_1,v_1\},\ldots,\{u_{2k+6},v_{2k+6}\}\}$, that is $z_c(u_i,v_i)\leq 3k$ for $1\leq i\leq 2k+6$.
Define a color set 
\[
C_Q=\cup_{p\in Q}\{c(T)\ |\ T\subseteq \mathcal{G}, p\subseteq T\}
\] and let $q=|C_Q|$. Then by Lemma \ref{pre}, we have the following  inequality $$q\leq 6k^2+18k.$$ Let $G$ be a rainbow subgraph of $\mathcal{G}$ with $c(n,k)-q$ edges and vertex set $[n]$ such that 
\[
\{c(e)\ |\ e\in E(G)\}\cap C_Q=\emptyset.
\]

\medskip
\textbf{Claim 1.}~$G$ is   $\mathcal{F}_k$-free.
\medskip

Otherwise, suppose that $G$ contains a copy of $\mathcal{F}_k$ denoted by $\mathcal{F}$ with center $u$. Since $|Q|\geq 2k+6$ and $|V(\mathcal{F})|=2k+1$, there exists $\{u_i,v_i\}\in Q$ such that $\{u_i,v_i\}\cap V(\mathcal{F})=\emptyset$. By the definition of $G$, $E(\mathcal{F})\cup \{u,u_i,v_i\}$ induces a rainbow copy of $\mathcal{F}_{k+1}$, contradicting the hypothesis. This completes the proof of Claim 1. 

Let  $w: G\times P \rightarrow R$ be a weighted function defined as in Section 2. And we denote $W_v=\sum_{p\in N(v)}w(v\cup p,p)$. Next we discuss two cases.

\medskip
\textbf{Case 1.~}$k$ is odd.
\medskip

Then we have
 \begin{equation}\label{e(G)-low}
    e(G)= c(n,k)-q\geq(n-2k)k(k-1)+2{k\choose 3}+2-6k^2-18k.
 \end{equation}
By Lemma \ref{weight}, we have $W_v\leq k(k-1)$ for all $v\in V(G)$. Let $W:=\{v\in V(G)\ |\ W_v=k(k-1)\}$. Then we have the following claim. 

\medskip
\textbf{Claim 2.~}$W\neq \emptyset$.
\medskip
 
Suppose to the contrary that $W_v<k(k-1)$ for all $v\in V(G)$. By the definition of weight function, we have $W_v\leq k(k-1)-2/3$ for all $v\in V(G)$. So 
\begin{align}\label{e(G)-upp}
e(G)=\sum_{v\in V(G)}W_v\leq nk(k-1)-2n/3.
\end{align}
Combining (\ref{e(G)-low}) and (\ref{e(G)-upp}), we may infer that
\[
n\leq \frac{5}{2}k^3+\frac{15}{2}k^2+26k-3,
\]
a contradiction. This completes the proof of Claim 2.

By Claim 2, we may choose $x\in W$ such that $W_x=k(k-1)$. Then by Lemma \ref{weight}, $N_G(x)$ is the disjoint union of two complete 2-graphs on $k$ vertices and every edge of $N_G(x)$ is in $A$-type.  Denote the two complete subgraphs by $R_1$ and $R_2$.

\medskip
\textbf{Claim 3.~} $|\{c(e)\ |\ x\in e\mbox{ and }e\subseteq {[n]\choose 3}\}|\leq 2{k\choose 2}+1$.
\medskip

Otherwise, suppose that $|\{c(e)\ |\ x\in e\mbox{ and }e\subseteq {[n]\choose 3}\}|\geq 2{k\choose 2}+2$. Then one can see that
\begin{align}\label{odd-2-color}
|\{c(p\cup \{x\})\ |\  p\notin E(R_1\cup R_2) \mbox{ and }p\in {[n]-x\choose 2}\}|\geq 2.    
\end{align}
We choose $T_1,T_2\in {[n]-\{x\}\choose 2}-({V(R_1)\choose 2}\cup {V(R_2)\choose 2})$ such that $c(T_1\cup \{x\})\neq c(T_2\cup \{x\})$, and for $i\in [2]$, $c(T_i\cup \{x\})\notin \{c(p\cup \{x\})\ |\ p\in E(R_1\cup R_2)\}$. For $i\in [2]$,  $F-V(T_i)$ has a matching $M_i$ of size $k-1$. Note that 
\[
|V(M_1\cup M_2)\cup T_1\cup T_2|\leq 2k+4.
\]
So by Lemma \ref{pre}, we may choose $T_0\in Q$ such that $x\notin T_0$ and $T_0\cap V(M_1\cup M_2\cup \{T_1,T_2\})=\emptyset.$ It follows that either $\{p\cup \{x\}\ |\ p\in M_1\cup \{T_1,T_0\}\}$ or $\{p\cup \{x\}\ |\ p\in M_2\cup \{T_2,T_0\}\}$ induces a rainbow copy of $\mathcal{F}_{k+1}$, a contradiction. 
This completes the proof of Claim 3.

Let $\mathcal{G}'$ be a rainbow edge-colored subgraph of  $\mathcal{G}-\{x\}$ obtained by deleting the triples colored by the colors from $\{c(e)\ |\ x\in e\mbox{ and }e\in {[n]\choose 3}\}$. 
By Claim 3, $e(\mathcal{G}')\geq f(n-1,k)+1$. Thus $\mathcal{G}'$ has a rainbow copy of $\mathcal{F}_k$ with core $y$. Choosing  $y'\in [n]-V(\mathcal{F}_k)-\{x\}$, in view of the preceding analysis, we know that $c(\{x,y,y'\})\notin \{c(e)\ |\ e\in \mathcal{F}_k\}$. It follows that $\{x,y,y'\}\cup \mathcal{F}_k$ a rainbow copy of $\mathcal{F}_{k+1}$, a contradiction. This completes the proof.

\medskip
\textbf{Case 2.~}$k\geq 4$ is even.
\medskip

Then we have 
\begin{align}\label{e(G)-low-even}
    e(G)\geq c(n,k)-q\geq \frac{1}{2}nk(2k-3)-\frac{1}{2}(2k^3-9k+6)-6k^2-18k.
\end{align}

By Lemma \ref{weight}, $W_v\leq k(k-3/2)$ for all $v\in V(G)$. By Theorem \ref{tang},  we may assume that $k\geq 4$.

\medskip
\textbf{Claim 4.~}There exists $v\in V(G)$ such that 
$W_v= k(k-3/2)$.
\medskip

Otherwise, we may assume that $W_v\leq  k(k-3/2)-1/2$ for all $v\in V(G)$. Thus we get 
\begin{align}\label{e(G)-up-even}
    e(G)\leq n k(k-3/2)-n/2.
\end{align}
Combining (\ref{e(G)-low-even}) and (\ref{e(G)-up-even}) , we may infer that $n\leq 2k^3+12k^2+27k+6$, a contradiction.
This completes the proof of Claim 4.

We choose $x\in V(G)$ with $W_x=k(k-3/2)$.
Suppose that  following inequality holds
\begin{align}\label{even-bound}
|\{c(e)\ |\ x\in e\mbox{ and }e\subseteq {[n]\choose 3}\}|\leq k(k-\frac{3}{2})+1,
\end{align}
Let $G'$ be a rainbow subgraph of $\mathcal{G}$ avoiding the colors appearing in  $\{c(e)\ |\ x\in e\mbox{ and }e\subseteq {[n]\choose 3}\}$. 
By (\ref{even-bound}), we have $e(G')>f(n-1,k)$. Hence $G'$ contains a copy  of $\mathcal{F}_k$ denoted by $\mathcal{M}$. Let  $y$ be the core of $\mathcal{M}$. We select a vertex 
 $x'\in [n]-V(\mathcal{M})-\{x\}$. Now $E(\mathcal{M})\cup \{\{x',x,y\}\}$ induces a rainbow copy of $\mathcal{F}_{k+1}$, a contradiction.  This contradiction concludes our proof. So next it is sufficient for us to show that (\ref{even-bound}) holds. 

Let $F:=N_{G}(x)$ denote the 2-graph with edge set $N_G(x)$ and vertex set $V(N_G(x))$. 
By Lemma  \ref{weight}, $F$ consists of two vertex-disjoint factor-critical graphs  or satisfies (a), (b), and (c).

We first prove the following claim, which  provides a method for constructing a rainbow $\mathcal{F}_{k+1}$ in $\mathcal{G}$ as described below.

\medskip
\textbf{Claim 5.~}Let $x_0 \in V(F)$ with $d_F(x_0) = k-1$, and let $\{y_0,w\} \subseteq [n]-(N_F(x_0)\cup \{x,x_0\})$ be a  2-element set. If $ c(\{x_0, w, x\}) \neq c(\{x_0, y_0, x\}) $, then $\mathcal{G}$ has a rainbow $F_{k+1} $ with core  $x_0$.
\medskip

  Let $N_F(x_0)=\{y_2,\ldots,y_{k}\}$. Our next step is to construct a rainbow copy  of $\mathcal{F}_{k+1}$ with the core $x_0$. According to Lemma \ref{weight}, for $2\leq i\leq k$,  $x_0y_i$ is $A$-type. Consequently, we know that 
 $d_G(\{x_0,y_i\})\geq 2k-1$.  For $2\leq i\leq k$, we can pick a set  $\{y_{i,1},...,y_{i,2k-1}\}\subseteq N_G(\{x_0,y_i\})$. Let $T_1$ be a 2-element subset chosen from the set $[n]-V(F)-\{x,w,y_0\}-(\bigcup_{i=2}^{k}\bigcup_{j=1}^{2k-1}\{y_{i,j}\})$. 
 Since $c(\{x_0,w,x\})\neq c(\{x_0,y_0,x\})$, we have either $c(T_1\cup\{x_0\})\neq c(\{x_0,y_0,x\})$ or $c(T_1\cup\{x_0\})\neq c(\{x_0,w,x\})$. Without loss of generality, assume that $c(T_1\cup\{x_0\})\neq c(\{x_0,y_0,x\})$. Since $G$ is a rainbow subgraph and $k\geq 4$, we can re-order the vertices $y_2,...,y_k$ in such a way that 
\[
c(T_1\cup\{x_0\})\notin \{c(\{x_0,y_{k-1},y_{k-1,j}\})\ |\ j\in[2k-1]\}
\] 
and
\[\{c(\{x_0,y_0,x\}),c(T_1\cup\{x_0\})\}\cap \{c(\{x_0,y_k,y_{k,j}\})\ |\ j\in[2k-1]\}=\emptyset.\]
 
Next, we aim to extend the sets $T_1\cup\{x_0\}$, $\{x_0,y_0,x\}$ to form a rainbow copy of $\mathcal{F}_{k+1}$ with the core $x_0$. Since $k\geq 4$,  we can find a vertex   $y_{2,1}'\in N_G(\{x_0,y_{2}\})-\{x, y_3,\ldots,y_k\}$, such that $\{\{T_1\cup\{x_0\},\{x_0,y_0,x\},\{x_0,y_{2,1}',y_2\}\}$ induces a rainbow copy $\mathcal{F}_3$. Now, assume that for some $i<k-1$, we have already constructed a collection $\{ \{T_1\cup\{x_0\},\{x_0,y_0,x\},\{x_0,y_{2,1}',y_2\},\ldots,\{x_0,y_{i,1}',y_i\} \}$ that induce a rainbow copy of $\mathcal{F}_{i+1}$, where $\{y_{2,1}',\ldots,y_{i,1}'\}\cap N_{F}(x_0)=\emptyset$.  Define $W_i:= \{\cup_{j=2}^i\{y_j,y_{j,1}'\}\}\cup \{y_0,x_0,x\}\cup V(T_1)$. Let's consider the step for $i+1$. Since $d_G(\{y_{i+1},x_0\})\geq 2k-1$ and $T_1\cap N_G(\{y_{i+1},x_0\})=\emptyset$, there exists  $y_{i+1,1}'\in N_G(\{y_{i+1},x_0\})-(W_i\cup N_{F}(x_0)$ 
 such that $\{\{T_1\cup\{x_0\},\{x_0,y_0,x\},\{x_0,y_{2,1}',y_2\},\ldots,\{x_0,y_{i+1,1}',y_{i+1}\} \}$ induces a rainbow copy of $\mathcal{F}_{i+2}$. Continuing the process until $i=k-1$, we can obtain a rainbow copy of $\mathcal{F}_{k}$ denoted by $\mathcal{M}$.  Recall that $d_G(\{y_k,x_0\})\geq 2k-1$ and $T_1\cap N_G(\{y_k,x_0\})=\emptyset$. As a result, there exists a vertex $y_{k,1}'$ in the set $ N_G(\{y_k,x_0\})-W_{k-1}$. Moreover, by the choice of $y_{k}$, 
$c(\{x_0,y_k,y_{k,1}'\})\notin \{c(\{x_0,y_0,x\}),c(T_1\cup \{x_0\})\}$. Hence $E(\mathcal{M})\cup \{x_0,y_k,y_{k,1}'\}$ induces a rainbow copy of $\mathcal{F}_{k+1}$, a contradiction.
This completes the proof of Claim 5. 

 Next by Lemma \ref{weight}, we proceed by discussing two subcases.
 
\medskip
\noindent\textbf{Subcase 2.1.~}$F=R_1\cup R_2$, where $R_1\cong K_{k-1}$ and $R_2$ is a factor-critical graph of order $k+1$  with degree sequence $(k-1,\ldots,k-1,k-2)$. 
\medskip


\medskip
\textbf{Claim 6.~} $\{c(g\cup \{x\})\ |\ g\in {[n]-\{x\}-V(R_1\cup R_2)\choose 2}\}$ contains exactly one color denoted by $c_0$ which belongs to the color set $C_Q$.
\medskip

 We will use a proof by contradiction. Firstly, suppose  that there exists $g\in {[n]-\{x\}-V(R_1\cup R_2)\choose 2}$ such that $c(g\cup \{x\})\in \{c(p\cup \{x\})\ |\ p\in E(R_1\cup R_2)\}$. Write $c(g\cup \{x\})=c(g'\cup \{x\})$, where $g'\in E(R_1\cup R_2)$. Since $R_1$ and $R_2$ are factor-critical graphs, $F-g'$ has a matching $M_1$ of size $k-1$. Let $h\in Q$ such that $x\notin h$ and $h\cap V(M_1\cup \{ g\})=\emptyset.$ It follows that $\{p\cup \{x\}\ |\ p\in M_1\cup \{g,h\}\}$ induces a rainbow copy of $\mathcal{F}_{k+1}$, a contradiction. Subsequently, suppose that  there are two distinct edges $g_1, g_2\in {[n]-\{x\}-V(R_1)-V(R_2)\choose 2}$ such that $c(g_1\cup \{x\})\neq c(g_2\cup \{x\})$ and neither of them belongs to $\{c(p\cup \{x\})\ |\ p\in E(R_1\cup R_2)\}$. 
Let $M_2$ be a matching  of size $k-1$ in $R_1\cup R_2$. Since $|Q|\geq 2k+6$, we may choose $h\in Q$ such that $x\notin h$ and $h\cap V(M_2\cup \{g_1,g_2\})=\emptyset.$ It follows that either $\{p\cup \{x\}\ |\ p\in M_2\cup \{g_1,h\}\}$ or $\{p\cup \{x\}\ |\ p\in M_2\cup \{g_2,h\}\}$ induces a rainbow copy of $\mathcal{F}_{k+1}$, a contradiction.  This completes the proof of Claim 6.

\medskip
\textbf{Claim 7.~}For any $f\in ([n]-V(R_1)-\{x\})\times V(R_1)$, we have $c(f\cup \{x\})\in \{c_0\}\cup \{c(p\cup \{x\})\ |\ p\in E(F)\}$;  and for any $ g\in ([n]-V(R_1\cup R_2)-\{x\})\times V(R_2)$, we have $c(g\cup \{x\})=c_0$. 
\medskip

 We will use a proof by contradiction.  
Firstly, suppose  that there exists   $f_1\in ([n]-V(R_1)-\{x\})\times V(R_1)$  such that  $c(f_1\cup \{x\})\notin \{c_0\}\cup \{c(p\cup \{x\})\ |\ p\in E(R_1\cup R_2)\}$. Since $R_1$ and $R_2$ are factor-critical graphs, $F-V(f_1)$ contains a matching $M_3$ of size $k-1$. Then we can select $h\in {[n]-V(R_1\cup R_2)-\{x\}\choose 2}$ such that $h\cap V(f_1)=\emptyset$. By Claim 6, we have $c(h\cup\{x\})=c_0$. Thus $\{p\cup \{x\}\ |\ p\in M_3\cup \{f_1,h\}\}$ induces a rainbow copy of $\mathcal{F}_{k+1}$, which is a contradiction. 
Next consider   that there exists  $g_1\in ([n]-V(R_1\cup R_2)-\{x\})\times V(R_2)$ such that $c(g_1\cup \{x\})\neq c_0$. If $c(g_1\cup \{x\})\notin \{c(p\cup \{x\})\ |\ p\in E(R_1\cup R_2)\}$, let $M_4$ be a matching of size $k-1$ of $F-V(g_1)$ and let $h'$ be a 2-subset of $[n]-(g_1\cup V(F)\cup \{x\})\cup V(M)$, then $\{p\cup \{x\}\ |\ p\in M_4\cup \{g_1,h'\}\}$ induces a rainbow copy of $\mathcal{F}_{k+1}$
 a contradiction again. So we may assume that  $c(g_1\cup \{x\})\in \{c(p\cup \{x\})\ |\ p\in E(R_1\cup R_2)\}$. Select an edge
 $g_2\in E(R_1\cup R_2)$ such that  $c(g_2\cup \{x\})=c(g_1\cup \{x\})$. Then, pick an element  $g_0\in {[n]-V(F)-V(g_1\cup g_2)-\{x\})\choose 2}$. By Claim 6, we have $c(g_0\cup \{x\})=c_0$.
If $F-V(g_1)-g_2$ has a matching $M$ of size $k-1$, then 
$\{e\cup \{x\}\ |\ e\in M\cup \{g_0,g_1\}\}$ induces a rainbow copy of $\mathcal{F}_{k+1}$, which is a contradiction. 
Consequently, we can assume that 
\begin{align}\label{F-ineq2}
  \nu(F-V(g_1)-g_2)<k-1.
\end{align}
 If $k\geq 6$, by  Lemma \ref{degree-critical}, $\nu(F-V(g_1)-g_2)=k-1$, which contradicts inequality (\ref{F-ineq2}). So we can conclude  $k=4$. 
Write $V(R_2):=\{y_0,y_1,\ldots,y_4\}$ such that $d_{R_2}(y_0)=2$, and $y_0y_1,y_0y_4\in E(R_2)$. If $g_1\cap \{y_1,y_4\}=\emptyset$ or $g_1\cap g_2\neq \emptyset$,  then  $F-V(g_1)-g_2$ has a matching of size $3$, which again contradicts inequality (\ref{F-ineq2}). Moreover, by symmetry, we may assume that $y_1\in g_1\cap V(R_2)$ and $g_1\cap g_2=\emptyset$.
Recall that for every edge $e\in E(R_2)$, $e$ is type $A$ in $G$, which means $d_{G}(\{y_1,y\})\geq 7$ for $y\in \{y_0,y_2,y_3\}$. Note that $\{y_2,y_3,y_0\}\cap (g_1\cup \{x\})=\emptyset$.  Let $w\in [n]-V(F)-\{x\}$. By Claim 7, we have $c(\{x,w,y_1\})=c_0\neq c(g_1\cup \{x\})$.  With the same discussion as Claim 6, 
$\{x\}\cup g_1$ can be extended into a rainbow copy with core $y_1$ denoted by $\mathcal{M}$ of $\mathcal{F}_4$, where $\{c(e)\ |\ e\in E(F_4)\}\subseteq \{c(e)\ |\ e\in E(\mathcal{M})\}$. By Lemma \ref{pre}, we may choose $h\in Q$ such that $h\cap V(\mathcal{M})=\emptyset$. Then $\{h\cup \{y_1\}\}\cup E(\mathcal{M})$ induces a rainbow copy of $\mathcal{F}_5$ in $\mathcal{G}$, a contradiction. 
This completes the proof of Claim 7. 

By Claims 6 and 7, it is sufficient for us to show that
 for any $l\in {V(R_2)\choose 2}\setminus E(R_2)$, $c(l\cup \{x\})=c_0$. 
  Suppose to the contrary that there exists  $l_0=\{x_0,y_0\}\in {V(R_2)\choose 2}\setminus E(R_2)$ such that  $c(\{x_0,y_0,x\})\neq c_0$. Note that $R_2$ is a factor-critical graph of order $k+1$  with degree sequence $(k-1,\ldots,k-1,k-2)$. Without loss of generality, suppose $d_{R_2}(x_0)= k-1.$ By Claim 7, for $w\in [n]-V(F)-\{x\}$, $c(\{x_0,w,x\})=c_0\neq c(\{x_0,y_0,x\})$. By Claim 5, we can find a rainbow copy of $\mathcal{F}_{k+1}$ with the core $x_0$, which contradicts the hypothesis. 
Thus we can derive that the inequality (\ref{even-bound}) holds.
 





\medskip
\noindent\textbf{Subcase 2.2.~}$F$ satisfies (a), (b) and (c). 
\medskip

Let \(U_k\) denote the set of vertices that serve as the centers of stars fulfilling condition (b) in Lemma \ref{weight}. Let \(V_k := V(F)\setminus(V(F_0)\cup U_k)\).  Let \(F_0\) be the factor-critical component that satisfies condition (a) in Lemma \ref{weight}. 
Since \(\nu(F)=k - 1\) and the number of vertices in \(F_0\) is at least \(k + 1\) (i.e., \(|F_0|\geq k + 1\)), we can infer that the cardinality of the set \(U_k\) is at most \(\frac{k}{2}-1\), that is, \(|U_k|\leq\frac{k}{2}-1\). 

 \medskip
\textbf{Claim 8.~}The set $\{c(g\cup \{x\})\ |\ g\in {[n]-\{x\}-V(F)\choose 2}\}$ contains exactly one color denoted by $c_1$, which belongs to the color set $C_Q$.  
\medskip

We will prove this statement by contradiction. First, suppose  that  there exists  $g\in {[n]-\{x\}-V(F)\choose 2}$ such that 
\[
c(g\cup \{x\})\in \{c(f\cup \{x\})\ |\ f\in E(F)\}.
\]
Our goal is to construct a rainbow copy   $\mathcal{F}_{k+1}$, which will lead to a contradiction.
Without loss of generality, assume that  $c(g\cup \{x\})=c(g'\cup \{x\})$, where $g'\in E(F)$. Since $F_0$ is factor-critical and $U_k\leq \frac{k}{2}-1$, $F-g'$ has a matching $M$ of size $k-1$. By Lemma \ref{pre}, we can then find  $h\in Q$ such that $x\notin h$ and $h\cap V(M\cup \{g\})=\emptyset.$ In this situation, the collection $\{p\cup \{x\}\ |\ p\in M\cup \{g,h\}\}$ forms a rainbow copy of $\mathcal{F}_{k+1}$, which contradicts our initial assumption. Next, consider the case when there are two distinct edges  $g_1, g_2\in {[n]-\{x\}-V(F)\choose 2}$ such that $c(g_1\cup \{x\})\neq c(g_2\cup \{x\})$ and neither $c(g_1\cup \{x\})$ nor $c(g_2\cup \{x\})$ belongs to the set $\{c(p\cup \{x\})\ |\ p\in E(F)\}$. 
We select a matching $M'$ of size $k-1$ in $F$. Subsequently, we can find an element  $h\in Q$ such that $x\notin h$ and $h\cap V(M'\cup \{g_1,g_2\})=\emptyset.$ As a consequence,  either the set $\{p\cup \{x\}\ |\ p\in M'\cup \{g_1,h\}\}$ or $\{p\cup \{x\}\ |\ p\in M'\cup \{g_2,h\}\}$ forms a rainbow copy of $\mathcal{F}_{k+1}$, a contradiction. This completes the proof of Claim 8.

Let $\mathcal{C}:=  \{c(p\cup \{x\})\ |\ p\in E(F)\}$ and let $\mathcal{C}':=\{c_1\}\cup  \{c(p\cup \{x\})\ |\ p\in E(F)\}.$

\medskip
\textbf{Claim 9.~}For any \(f\in \binom{[n]\setminus (U_k\cup\{x\})}{2}\) satisfying  \(|f\cap V(F_0)| = 1\), we have \(c(f\cup \{x\})=c_1\). 
\medskip

Suppose  that there exists $g_1\in {[n]-U_k-\{x\}\choose 2}$ such that $|g_1\cap V(F_0)|=1$ and  $c(g_1\cup \{x\})\neq c_1$.  We choose $h\in {[n]-V(F)-\{x\}\choose 2}$ such that $h\cap V(g_1)=\emptyset$.  By Claim 8, we have $c(h\cup\{x\})=c_1$. If $c(g_1\cup \{x\})\notin \mathcal{C}$, we proceed as follows. Let $M_3$ be a matching $M_3$ of size $k-1$ in $F-V(g_1)$. Then, the collection $\{p\cup \{x\}\ |\ p\in M_3\cup \{g_1,h\}\}$ induces a rainbow copy of $\mathcal{F}_{k+1}$, a contradiction.  
Hence we may infer that  
\begin{align}
c(g_1\cup \{x\})\in \mathcal{C}.
\end{align}
Let $g_2\in E(F)$ such that $c(g_2\cup \{x\})=c(g_1\cup \{x\})$.
If $F-V(g_1)-g_2$ contains a matching $M_4$ of size $k-1$, $\{p\cup \{x\}\ |\ p\in M_4\cup \{g_1,h\}\}$ induces rainbow copy of $\mathcal{F}_{k+1}$, a contradiction.  
Thus we have 
\begin{align}\label{c11-eq2}
\nu(F-V(g_1)-g_2)<k-1.
\end{align}

Then by Lemma \ref{degree-critical}, we may infer that $k= 4$. 
Let $g_1\cap V(F_0)=\{y_0\}$. Consider $d_F(y_0)=3$. 
By the proof of Claim 5, $\{x\}\cup g_1$ may be extended into a rainbow copy denoted by $\mathcal{M}$ of $\mathcal{F}_4$ such that for any $e\in E(\mathcal{M})$, $c(e)\in \{c(f)\ |\ f\in E(G)\}$. 
Now we select $h\in Q$ such that $h\cap V(\mathcal{M})=\emptyset$.  It follows that $(\{\{y_0\}\cup h\})\cup E(M)$ induces a rainbow copy of $\mathcal{F}_5$, a contradiction. 
Next we may assume $d_F(y_0)=2$. Recall that $g_2\in E(F)$ such that 
$c(g_2\cup \{x\})= c(g_1\cup \{x\})$. By Lemma \ref{cycle}, $F-V(g_1)$ contains a Hamilton cycle. So 
$F-V(g_1)-g_2$ contains a matching $M_5$ of size $3$, which contradicts the inequality (\ref{c11-eq2}). This completes the proof of Claim 9. 


 \medskip
\textbf{Claim 10.~}For any \(f\in \binom{[n]\setminus (U_k\cup\{x\}\cup V(F_0))}{2}\) satisfying  \(|f\cap V_k| = 1\), we have \(c(f\cup \{x\})=c_1\).
\medskip

We proceed by contradiction. 
Suppose  to the contrary that there exists  $f_1\in \binom{[n]\setminus \big(U_k\cup\{x\}\cup V(F_0)\big)}{2}$ such that $c(f_1\cup \{x\})\neq c_1$.
 We choose $h\in {[n]-V(F)-\{x\}\choose 2}$ such that $h\cap V(f_1)=\emptyset$. It follows that $c(h\cup\{x\})=c_1$ by Claim 8.  If  $c(f_1\cup \{x\})\in \mathcal{C}$, let $f_2\in E(F)$ such that $c(f_1\cup \{x\})=c(f_2\cup \{x\})$; otherwise let $f_2=\emptyset$. Since $F_0$ is a factor-critical graph, we know that $F-V(f_1)-f_2$ contains a matching $M_2$ of size $k-1$. 
Thus  $\{p\cup \{x\}\ |\ p\in M_2\cup \{f_1,h\}\}$ induce a rainbow copy of $\mathcal{F}_{k+1}$, a contradiction again. This completes the proof of Claim 10.

 \medskip
\textbf{Claim 11.~}For any $ T\in {V_k\choose 2}$, $c(T\cup \{x\})\in \mathcal{C}'$. 
\medskip

  Suppose to the contrary that there exists $T_1\in {V_k\choose 2}$ such that $c(T_1\cup \{x\})\notin \mathcal{C}'$. Since $d_F(w)= 3$ for  $w\in U_k$ and $k\geq 4$, $F-V(T_1)$ has a matching $M$ of size $k-1$. We may choose $T_2\in {[n]-(\{x\}\cup V(F)\cup T_1)\choose 2}$. It follows that $c(T_2\cup\{x\})=c_1$ by Claim 8. So $\{p\cup \{x\}\ |\ p\in M\cup \{T_2,T_1\}\}$ induces a rainbow copy of $\mathcal{F}_{k+1}$, a contradiction.  This completes the proof of Claim 11. 

Let $z\in V(F_0)$ such that $d_{F_0}(z)=k-2$. 

\medskip
\textbf{Claim 12.~} For any $T\in \left((U_k\times V(F_0-z)) \cup {V(F_0)\choose 2}\right)\setminus E(F)$, $c(T\cup\{x\})\in \mathcal{C}'$.
\medskip

 Suppose to the contrary that there exists $T\in \left((U_k\times V(F_0-z)) \cup {V(F_0)\choose 2}\right)\setminus E(F)$ such that $c(T\cup \{x\})\notin \mathcal{C}'$.
 Write $T=\{y_1,y_2\}$. Without loss of generality, let $y_1$ be a vertex in $F_0$ with degree $k-1$. We choose $w\in [n]-(V(F)\cup \{x\})$. By Claim 9, we have $c(\{x,w,y_1\})=c_1$. So we have 
  $c(\{x,w,y_1\})\neq c(\{x,y_1,y_2\})$. Recall that $\{x,y_1\},\{y_1,y_2\}\notin E(F)$. Then by Claim 5, 
 we can find a rainbow copy of $\mathcal{F}_{k+1}$  in $\mathcal{G}$ with core $y_1$, a contradiction. This completes the proof of Claim 12. 

\medskip
\textbf{Claim 13.~}For any $T\in {[n]-(V(F_0-z)\cup \{x\})\choose 2}\setminus E(F)$, $c(T\cup \{x\})\in \mathcal{C}'$.  
\medskip

Suppose to the contrary that there exists $T\in {[n]-(V(F_0)\cup \{x\})\choose 2}\setminus E(F)$ such that $c(T\cup \{x\})\notin \mathcal{C}'$.  Write $T=\{x_0,y_0\}$. By Claims 8, 9, 10 and 11, we have $T\cap U_k\neq \emptyset$.  Without loss generality, suppose that $x_0\in U_k$. Then  $d_{F}(x_0)= k-1.$  
Let $w\in V(F_0)-N_F(x_0)$ such that $d_{F_0}(w)=k-1$. By Claim 12, we have $c(\{w,x,x_0\})\in \{c(p\cup\{x\})\ |\ p\in E(F)\}\cup \{c_1\}$. Thus it follows that $c(\{w,x,x_0\})\neq c( \{x,x_0,y_0\})$. By Claim 5, we can find a rainbow copy of $\mathcal{F}_{k+1}$ with core $x_0$, a contradiction. This completes the proof of Claim 13.

By Claims 8,9,10,11,12 and 13, for any $e\in {[n]-x\choose 2}$, $c(e\cup \{x\})\in \mathcal{C'}$. Thus we can derive that the inequality (\ref{even-bound}) holds. This complete the proof of Theorem \ref{main}. \qed

\end{document}